\documentclass{ article}

\begin{document}
\vspace*{.5cm}
\begin{center}
{\Large{\bf   Slant Riemannian maps from almost Hermitian manifolds}}\\
\vspace{.5cm}
 { Bayram \d{S}ahin} \\
\end{center}

\vspace{.5cm}
\begin{center}
{\it Inonu University, Department of Mathematics, 44280,
Malatya-Turkey. E-mail:bayram.sahin@inonu.edu.tr}
\end{center}
\vspace{.5cm}

\noindent {\bf Abstract.} {\small  As a generalization of holomorphic submersions, anti-invariant submersions and slant submersions, we introduce slant Riemannian maps from almost Hermitian manifolds to Riemannian manifolds. We give examples, obtain the existence conditions of slant Riemannian maps and investigate harmonicity of such maps. We also obtain necessary and sufficient conditions for slant Riemannian maps to be totally geodesic and give a decomposition theorem for the total manifold.}\\



\pagestyle{myheadings}

\section*{1.~Introduction}

  \setcounter{equation}{0}
\renewcommand{\theequation}{1.\arabic{equation}}
\markboth{Slant Riemannian maps}{Slant Riemannian
maps}{\thispagestyle{plain}}

Differentiable maps between Riemannian manifolds are important in differential geometry. There are certain types of differentiable maps between Riemannian manifolds whose existence influence the geometry of the source manifolds and the target manifolds. Differentiable maps between Riemannian manifolds are also useful to compare geometric structures defined on both manifolds. Basic maps in this manner
are isometric immersions between Riemannian manifolds. Such maps are
characterized by their Jacobian matrices and the
induced metric which is symmetric positive definite  bilinear form.
The theory of isometric immersions is an active research area  and
it plays an important role in the development of modern differential
geometry. The other basic maps for comparing geometric structures
defined on Riemannian manifolds are Riemannian submersions and they
were studied by O'Neill \cite{O'Neill} and Gray \cite{Gray}. The theory of Riemannian submersions is also a very active research field, for recent developments in this area see:\cite{Falcitelli}.\\

In 1992, Fischer  introduced Riemannian maps between Riemannian
manifolds in \cite{Fischer} as a generalization of the notions of
isometric immersions and Riemannian submersions. Let $F:(M_1,
g_1)\longrightarrow (M_2, g_2)$  be a smooth map between Riemannian
manifolds such that $0<rank F<min\{ m, n\}$, where $dimM_1=m$ and
$dimM_2=n$. Then we denote the kernel space of $F_*$ by $kerF_*$ and
consider the orthogonal complementary space
$\mathcal{H}=(kerF_*)^\perp$ to $kerF_*$ in $TM_1$. Thus the tangent
bundle of $M_1$ has the following decomposition

$$TM_1=kerF_* \oplus\mathcal{H}.$$

We denote the range of $F_*$ by $rangeF_*$ and consider the
orthogonal complementary space $(rangeF_*)^\perp$ to $rangeF_*$ in
the tangent bundle $TM_2$ of $M_2$. Since $rankF<min\{ m, n\}$, we
always have $(rangeF_*)^\perp$. Thus the tangent bundle $TM_2$ of
$M_2$ has the following decomposition
$$TM_2=(rangeF_*)\oplus (rangeF_*)^\perp.$$
 Now, a smooth
map $F:(M^{^m}_1,g_1)\longrightarrow (M^{^n}_2, g_2)$ is called
Riemannian map at $p_1 \in M_1$ if the horizontal restriction
$F^{^h}_{*p_1}: (ker F_{*p_1})^\perp \longrightarrow (range
F_{*p_1})$  is a linear isometry between the inner product spaces
$((ker F_{*p_1})^\perp, g_1(p_1)\mid_{(ker F_{*p_1})^\perp})$ and
$(range F_{*p_1}, g_2(p_2)\mid_{(range F_{*p_1})})$, $p_2=F(p_1)$.
Therefore Fischer stated in \cite{Fischer} that a Riemannian map
is a map which is as isometric as it can be. In another words,
$F_*$ satisfies the equation
\begin{equation}
g_2(F_*X, F_*Y)=g_1(X, Y)\label{eq:1.1}
\end{equation}
 for $X, Y$ vector
fields tangent to $\mathcal{H}$. It follows that isometric
immersions and Riemannian submersions are particular Riemannian maps
with $kerF_*=\{ 0 \}$ and $(rangeF_*)^\perp=\{ 0 \}$. It is known
that a Riemannian map is a subimmersion \cite{Fischer} and this fact
implies that  the rank of the linear map $F_{*p}:T_pM_1 \longrightarrow T_{F(p)}M_2$ is constant for $p$ in each connected component of $M_1$, \cite{Marsden} and \cite{Fischer}. It is also important to note that Riemannian maps satisfy the eikonal equation which is a bridge between geometric optics and physical optics. For Riemannian maps and their applications in spacetime geometry, see: \cite{Garcia-Rio-Kupeli}.\\

Let $\bar{M}$ be a K\"{a}hler manifold with complex structure $J$
and $M$ a Riemannian manifold isometrically immersed in $\bar{M}$. A
submanifold $M$ is called holomorphic (complex) if $J (T_p M)\subset
T_p M$, for every $p \in M$, where $T_p M$  denotes the tangent
space to $M$ at the point $p $. $M$ is called totally real if $J(T_p
M) \subset T_p M^{\perp}$ for every $p \in M,$ where $T_p M^{\perp}$
denotes the normal space to $M$ at the point $p$. On the other hand,
a submanifold $M$ is called slant if for all
    non-zero vector $X$ tangent to $M$ the angle $\theta (X)$
    between $J X$ and $T_p M$ is a constant, i.e, it does not depend
    on the choice of $p \in M$ and $X \in T_p M$ \cite{Chen}. Holomorphic submanifolds and totally real submanifolds are slant
submanifolds with $\theta=0$ and $\theta=\frac{\pi}{2}$, respectively. A slant submanifold is called proper if it is neither holomorphic nor totaly real. \\

Riemannian  submersions between Riemannian manifolds  equipped with
differentiable structures were studied by Watson in \cite{Watson}.
As an analogue of holomorphic submanifolds, Watson defined almost
Hermitian submersions between almost Hermitian manifolds as follows:
Let $M$ be a complex $m-$dimensional almost Hermitian manifold with
Hermitian metric $g_M$ and almost complex structure $J_M$ and $N$ be
a complex $n-$dimensional almost Hermitian manifold with Hermitian
metric $g_N$ and almost complex structure $J_N$. A Riemannian
submersion $F:M\longrightarrow N$ is called an almost Hermitian
submersion if $F$ is an almost complex mapping, i.e.,
$F_*J_M=J_NF_*$. The main result of this notion is that the vertical
and horizontal distributions are $J_M-$ invariant. Watson also
showed that the base manifold and each fiber have the same kind of
structure as the total space, in most cases \cite{Watson} and
\cite{Falcitelli}. Since then almost Hermitian submersions have been
extended to the almost contact manifolds
\cite{Domingo}, \cite{Ianus2}, locally conformal K\"{a}hler
manifolds \cite{lck} and
quaternion K\"{a}hler manifolds \cite{Ianus}.\\

In \cite{Sahin}, we introduced anti-invariant Riemannian submersions
from almost Hermitian manifolds onto Riemannian manifolds  as
follows. Let $M$ be a complex $m-$ dimensional almost Hermitian
manifold with Hermitian metric $g_{_M}$ and almost complex structure
$J$ and $N$ be a Riemannian manifold with Riemannian metric
$g_{_N}$. Suppose that there exists a Riemannian submersion $F:M
\longrightarrow N$ such that the integral manifold of the
distribution $ker F_*$ is anti-invariant with respect to $J$, i.e.,
$J(ker F_*)\subseteq (ker F_*)^{\perp}$. Then we say that $F$ is an
anti-invariant Riemannian submersion. As a generalization of almost
Hermitian submersions and anti-invariant Riemannian submersions,
recently, we also introduced the notion of slant submersions from
almost Hermitian manifolds onto arbitrary Riemannian manifolds
\cite{Sahin1} as follows: Let $F$ be a Riemannian submersion from an
almost Hermitian manifold $(M_1,g_1,J_1)$  onto a Riemannian
manifold $(M_2,g_2)$. If for any non-zero vector $X \in \Gamma(ker
F_*)$, the angle $\theta(X)$ between $JX$ and the space $ker F_*$ is
a constant, i.e. it is independent of the choice of the point $p \in
M_1$ and choice of the tangent vector $X$ in $ker F_*$, then we say
that $F$ is a slant submersion. In this case, the angle $\theta$ is
called the slant
angle of the slant submersion. \\

In \cite{Sahin4}, as a generalization of almost Hermitian submersions, anti-invariant Riemannian submersions and slant submersions, we defined semi-invariant Riemannian maps from almost Hermitian manifolds and investigated the geometry of the total manifold and the base manifold by using the existence of such maps.\\

In this paper, as another generalization of Hermitian submersions,
anti-invariant submersions and slant submersions, we define and
study slant Riemannian maps from almost Hermitian manifolds to
Riemannian manifolds. In section 2, we recall basic facts for
Riemannian maps and almost Hermitian manifolds. In section 3, we
define slant Riemannian maps and give many examples. We also obtain
a characterization of such maps and investigate the harmonicity of
slant Riemannian maps. Then we give necessary and sufficient
conditions for slant Riemannian maps to be totally geodesic.
Finally, in section 4, we obtain a decomposition theorem for the
total manifold by using slant Riemannian maps.

\section*{\bf 2.~Riemannian maps}
\setcounter{equation}{0}
\renewcommand{\theequation}{2.\arabic{equation}}

In this section, we develop fundamental formulas for Riþemannian
maps similar to the Gauss-Weingarten formulas of isometric
immersions and O'Neill's formulas of Riemannian submersions. We
also recall useful results which are related to the second
fundamental form and the tension field of Riemannian maps. Let
$(M, g_{_M})$ and $(N, g_{_N})$ be Riemannian manifolds and
suppose that $F: M\longrightarrow N$ is a smooth map between them.
Then the differential $F_*$ of $F$ can be viewed a section of the
bundle $Hom(TM, F^{-1}TN)\longrightarrow M,$ where $F^{-1}TN$ is
the pullback bundle which has fibres $(F^{-1}TN)_p=T_{F(p)} N, p
\in M.$ $Hom(TM, F^{-1}TN)$ has a connection $\nabla$ induced from
the Levi-Civita connection $\nabla^M$ and the pullback connection.
Then the second fundamental form of $F$ is given by
\begin{equation}
(\nabla {F}_*)(X, Y)=\nabla^{F}_X {F}_*(Y)-{F}_*(\nabla^M_X Y)
\label{eq:2.1}
\end{equation}
for $X, Y \in \Gamma(TM).$ It is known that the second fundamental
form is symmetric \cite{Baird-Wood}.  First note that in
\cite{Sahin5} we showed that the second fundamental form $(\nabla
F_*)(X,Y)$, $\forall X, Y \in \Gamma((ker F_*)^\perp)$, of a
Riemannian map has no components in $range
F_*$. More precisely we have the following.\\

\noindent{\bf Lemma~2.1.~}{\it Let $F$ be a Riemannian map from a
Riemannian manifold $(M_1,g_1)$ to a Riemannian manifold
$(M_2,g_2)$. Then}
$$
g_2((\nabla F_*)(X,Y),F_*(Z))=0, \forall X,Y,Z \in \Gamma((ker
F_*)^\perp).
$$

As a result of Lemma 2.1, we have
\begin{equation}
(\nabla F_*)(X,Y) \in \Gamma((range F_*)^\perp), \forall X,Y \in
\Gamma((ker F_*)^\perp). \label{eq:2.2}
\end{equation}

For the tension field of a Riemannian map between Riemannian manifolds, we have the following.\\

\noindent{\bf Lemma~2.2.~}\cite{Sahin3}{\it Let $F: (M, g_{_M})
\longrightarrow (N,g_{_N})$ be a Riemannian map between Riemannian
manifolds. Then the tension field $\tau$ of $F$ is}
\begin{equation}
\tau=-m_1F_*(H)+m_2H_2, \label{eq:2.5}
\end{equation}
{\it where $m_1=dim(ker F_*), m_2=rank F$, $H$ and
$H_2$ are the mean curvature vector fields of the distribution $ker F_*$ and $range F_*$, respectively.}\\

Let $F$ be a Riemannian map from a Riemannian manifold $(M_1,g_1)$
to a Riemannian manifold $(M_2,g_2)$. Then we define $\mathcal{T}$
and $\mathcal{A}$ as

\begin{equation}
\mathcal{A}_E
F=\mathcal{H}\nabla_{\mathcal{H}E}\mathcal{V}F+\mathcal{V}\nabla_{\mathcal{H}E}\mathcal{H}F
\label{eq:2.3}
\end{equation}
\begin{equation}
\mathcal{T}_E
F=\mathcal{H}\nabla_{\mathcal{V}E}\mathcal{V}F+\mathcal{V}\nabla_{\mathcal{V}E}\mathcal{H}F,
\label{eq:2.4}
\end{equation}
for vector fields $E, F$ on $M_1$, where $\nabla$ is the
Levi-Civita connection of $g_1$. In fact, one can see that these
tensor fields are O'Neill's tensor fields which were defined for
Riemannian submersions. For any $E \in \Gamma(TM_1)$,
$\mathcal{T}_E$ and $\mathcal{A}_E$ are skew-symmetric operators
on $(\Gamma(TM_1),g)$ reversing the horizontal and the vertical
distributions. It is also easy to see that $\mathcal{T}$ is
vertical, $\mathcal{T}_E=\mathcal{T}_{\mathcal{V}E}$ and
$\mathcal{A}$ is horizontal,
$\mathcal{A}=\mathcal{A}_{\mathcal{H}E}$. We note that the tensor
field $\mathcal{T}$  satisfies
\begin{equation}
\mathcal{T}_UW=\mathcal{T}_WU, \forall U,W \in \Gamma(ker F_*).\label{eq:2.5}\\
\end{equation}
On the other hand, from (\ref{eq:2.3}) and (\ref{eq:2.4}) we have
\begin{eqnarray}
\nabla_VW&=&\mathcal{T}_VW+\hat{\nabla}_V W \label{eq:2.7}\\
\nabla_VX&=&\mathcal{H}\nabla_VX+\mathcal{T}_VX \label{eq:2.8}\\
\nabla_XV&=&\mathcal{A}_XV+\mathcal{V}\nabla_XV \label{eq:2.9}\\
\nabla_XY&=&\mathcal{H}\nabla_XY+\mathcal{A}_XY \label{eq:2.10}
\end{eqnarray}
for $X,Y \in \Gamma((ker F_*)^{\perp})$ and $V,W \in \Gamma(ker
F_*)$, where $\hat{\nabla}_VW=\mathcal{V}\nabla_VW$. \\

From now on, for simplicity, we denote by $\nabla^2$ both the
Levi-Civita connection of $(M_2, g_2)$ and its pullback along $F$.
 Then according to \cite{Nore}, for any vector field $X$ on $M_1$ and any section $V$ of $(range F_*)^\perp$, where $(range F_*)^\perp$ is the
 subbundle of $F^{-1}(TM_2)$ with fiber $(F_*(T_pM))^\perp$-orthogonal complement of $F_*(T_pM)$ for $g_2$ over $p$, we have
 $\nabla^{^{F \perp}}_XV$ which is the orthogonal projection of $\nabla^2_XV$ on $(F_*(TM))^\perp$. In \cite{Nore}, the author also showed that $\nabla^{^{F \perp}}$ is a linear connection on $(F_*(TM))^\perp$ such that $\nabla^{^{F \perp}}g_2=0$. We now  define $\mathcal{S}_{V}$ as
\begin{equation}
\nabla^2_{_{F_*X}}V=-\mathcal{S}_{_V}F_*X+\nabla^{^{F
\perp}}_{_{X}}V, \label{eq:2.11}
\end{equation}
where $\mathcal{S}_{_V}F_*X$ is the tangential component (a vector
field along $F$) of $\nabla^2_{_{F_*X}}V$. It is easy to see that
$\mathcal{S}_V F_*X$ is bilinear in $V$ and $F_*X$ and
$\mathcal{S}_V F_*X$ at $p$ depends only on $V_p$ and $F_{*p}X_p$.
By direct computations, we obtain
\begin{equation}
g_2(\mathcal{S}_{_V} F_*X,F_*Y)=g_2(V, (\nabla F_*)(X,Y)),
\label{eq:2.12}
\end{equation}
for $X, Y \in \Gamma((ker F_*)^\perp)$ and $V \in \Gamma((range
F_*)^\perp)$. Since $(\nabla F_*)$ is symmetric, it follows that
$\mathcal{S}_{_V}$ is a symmetric linear transformation of $range F_*$.\\

A $2k$-dimensional Riemannian manifold $(\bar{M}, \bar{g},
\bar{J})$ is called an almost Hermitian manifold if there exists a
tensor filed  $\bar{J}$ of type (1,1) on $\bar{M}$ such that
$\bar{J}^2 = - I$ and
\begin{equation}
\bar{g}(X,Y)=\bar{g}(\bar{J}X,\bar{J}Y) \label{eq:2.10}, \forall
X,Y \in \Gamma(T\bar{M}),
\end{equation}
 where $I$ denotes the identity transformation of
$T_{p} \bar{M}$. Consider an almost Hermitian manifold $(\bar{M},
\bar{J}, \bar{g})$  and denote by $\bar{\nabla}$ the Levi-Civita
connection on $\bar{M}$ with respect to $\bar{g}.$ Then $\bar{M}$
is called a K\"{a}hler manifold if $\bar{J}$  is
parallel with respect to $\bar{\nabla}$, i.e,
\begin{equation}
(\bar{\nabla}_X \bar{J})Y=0 \label{eq:2.11}
\end{equation}
for $X, Y \in \Gamma(T\bar{M}) \cite{Yano-Kon}.$

\section*{3.~Slant Riemannian maps}
  \setcounter{equation}{0}
\renewcommand{\theequation}{3.\arabic{equation}}

In this section, as a generalization of almost Hermitian submersions, slant submersions and
anti-invariant Riemannian submersions,  we introduce slant
Riemannian maps from an almost Hermitian manifold to a Riemannian
manifold. We first focus on the existence of such maps by giving
some examples. Then we investigate the effect of slant Riemannian maps on the
geometry of the total manifold, the base manifold and themselves.
More precisely, we investigate the geometry of leaves of
distributions on the total manifold arisen from such maps. We also
obtain necessary and sufficient conditions for slant Riemannian maps
to be harmonic and totally geodesic. We first present the following
definition.\\

\noindent{\bf Definition~3.1.~}{\it Let $F$ be a Riemannian map
from an almost Hermitian manifold $(M_1,g_1,J_1)$  to a Riemannian
manifold $(M_2,g_2)$. If for any non-zero vector $X \in \Gamma(ker
F_*)$, the angle $\theta(X)$ between $JX$ and the space $ker F_*$
is a constant, i.e. it is independent of the choice of the point
$p \in M_1$ and choice of the tangent vector $X$ in $ker F_*$,
then we say that $F$ is a slant Riemannian map. In this case, the angle $\theta$ is called the slant angle of the slant Riemannian map.}\\

Since $F$ is a subimmersion, it follows that the rank of $F$ is
constant on $M_1$, then the rank theorem for functions implies that
$ker F_*$ is an integrable subbundle of $TM_1$, (\cite{Marsden},
page:205).  Thus it follows from above definition that the leaves of
the distribution $ker F_*$ of a slant Riemannian map are slant
submanifolds of $M_1$, for
slant submanifolds, see: \cite{CB}.\\

We first give some examples of slant Riemannian maps.\\

\noindent{\bf Example~1.~}Every Hermitian submersion from an
almost Hermitian manifold onto an almost Hermitian manifold is a
slant
Riemannian map with $\theta=0$ and $(range F_*)^\perp=\{0\}$.\\

\noindent{\bf Example~2.~}Every anti-invariant Riemannian
submersion from an almost Hermitian manifold onto a Riemannian
manifold is a slant Riemannian map with $\theta=\frac{\pi}{2}$ and $(range F_*)^\perp=\{0\}$.\\

\noindent{\bf Example~3.~} Every proper slant submersion with the
slant angle $\theta$ is a slant Riemannian map with $(range
F_*)^\perp=\{0\}$.\\

We now denote the Euclidean $2m-$ space with the
standard metric by  $R^{2m}$. An almost complex structure $J$ on $R^{2m}$ is said
to be compatible if $(R^{2m},J)$ is complex analytically isometric
to the complex number space $C^m$ with the standard flat
K\"{a}hlerian metric. Then the compatible almost complex
structure  $J$ on $R^{2m}$ defined by
$$J(a^1,...,a^{2m})=(-a^{m+1},-a^{m+2},...,-a^{2m},a^1,a^2,...,a^m).$$

A slant Riemannian map is said to be proper if it is not a
submersion.
Here is an example of proper slant Riemannian maps.\\

 \noindent{\bf Example~4.~}Consider the following Riemannian
map given by
$$
\begin{array}{cccc}
  F: & R^4             & \longrightarrow & R^4\\
     & (x_1,x_2,x_3,x_4) &             & (0,\frac{x_2\,\sin\, \alpha +x_3+x_4\cos \alpha}{\sqrt{2}},0,x_2\,\cos\, \alpha-x_4\sin\, \alpha).
\end{array}
$$
Then for any $0< \alpha <\frac{\pi}{2}$, $F$ is a slant Riemannian
map with respect to the compatible almost complex
structure  $J$ on $R^{4}$ with slant angle $\frac{\pi}{4}$.\\

Let $F$ be a  Riemannian map from a K\"{a}hler manifold
$(M_1,g_1,J)$ to a Riemannian manifold $(M_2,g_2)$. Then for $X \in
\Gamma(ker F_*)$, we write
\begin{equation}
JX=\phi X+\omega X, \label{eq:3.1}
\end{equation}
where $\phi X$ and $ \omega X$ are vertical and horizontal parts
of $JX$. Also for $V \in \Gamma((ker F_*)^\perp)$, we have
\begin{equation}
JZ=\mathcal{B}Z+\mathcal{C}Z, \label{eq:3.2}
\end{equation}
where $\mathcal{B}Z$ and $\mathcal{C}Z$ are vertical and
horizontal components of $JZ$. Using (\ref{eq:2.7}),
(\ref{eq:2.8}), (\ref{eq:3.1}) and (\ref{eq:3.3}) we obtain
\begin{eqnarray}
(\nabla_X \omega)Y&=&\mathcal{C}\mathcal{T}_XY-\mathcal{T}_X \phi
Y
\label{eq:3.3}\\
(\nabla_X \phi)Y&=&\mathcal{B}\mathcal{T}_XY-\mathcal{T}_X \omega
Y, \label{eq:3.4}
\end{eqnarray}
where $\nabla$ is the Levi-Civita connection on $M_1$ and
\begin{eqnarray}
(\nabla_X \omega)Y&=&\mathcal{H}\nabla_X \omega Y-\omega
\hat{\nabla}_X Y\nonumber\\
(\nabla_X \phi)Y&=& \hat{\nabla}_X \phi Y -\phi
\hat{\nabla}_XY\nonumber
\end{eqnarray}
for $X, Y \in \Gamma(ker F_*)$. Let $F$ be a slant Riemannian map
from an almost Hermitian manifold $(M_1,g_1,J_1)$ to a Riemannian
manifold $(M_2,g_2)$, then we say that $\omega$ is parallel with
respect to the Levi-Civita connection $\nabla$ on $ker F_*$ if its
covariant derivative with respect to $\nabla $ vanishes, i.e., we
have
$$(\nabla_X \omega)Y=\nabla_X \omega Y-\omega(\nabla_X Y) =0$$
for $X, Y \in \Gamma(ker F_*)$. Let $F$ be a slant Riemannian map
from a complex $m$-dimensional Hermitian manifold $(M,g_1,J)$ to a
Riemannian manifold $(N,g_2)$. Then, $\omega(ker F_*)$ is a subspace
of $(ker F_*)^\perp$.  Thus it follows that $kerF_{*p} \oplus
\omega(ker F_{*p})$ is invariant with respect to $J$. Then for every
$p \in M$, there exists an invariant subspace $\mu_p$ of $(ker
F_{*p})^\perp$ such that
$$T_pM=kerF_{*p}\oplus \omega(ker F_{*p})\oplus \mu_p.$$

The proof of the following result is exactly the same with slant
immersions (see \cite{Chen} or \cite{Carriazo} and \cite{Carriazo2}  for Sasakian case),
therefore we omit its proof.\\

\noindent{\bf Theorem~3.1.~}{\it Let $F$ be a Riemannian map from an
almost Hermitian manifold $(M_1,g_1,J)$ to a Riemannian manifold
$(M_2,g_2)$. Then $F$ is a slant Riemannian map if and only if there
exists a constant $\lambda \in [-1,0]$ such that}
$$\phi^2X=\lambda X$$
{\it for $X \in \Gamma(ker F_*)$. If $F$ is a slant Riemannian map,
then
$\lambda=-\cos^2 \theta$.}\\

By using above theorem, it is easy to see that
\begin{eqnarray}
g_1(\phi X, \phi Y)&=&\cos^2 \theta g_1(X,Y) \label{eq:3.5}\\
g_1(\omega X, \omega Y)&=&\sin^2 \theta g_1(X,Y) \label{eq:3.6}
\end{eqnarray}
for any $X, Y \in \Gamma(ker F_*)$. Also by using (\ref{eq:3.5})
we can easily conclude that
$$\{e_1,\sec \theta \phi e_1, e_2,\sec \theta \phi e_2,...,e_n,\sec \theta \phi
e_n\}$$ is an orthonormal frame for $\Gamma(ker F_*)$. On the
other hand,  by using (\ref{eq:3.6}) one can see that
$$\{\csc
\theta \omega e_1, \csc \theta \omega e_2,...,\csc \theta \omega
e_n\}$$ is an orthonormal frame for $\Gamma(\omega(ker F_*))$. As in
slant immersions, we call  the frame
$$\{e_1,\sec \theta \phi e_1, e_2,\sec \theta \phi e_2,...,e_n,\sec \theta \phi
e_n,\csc \theta \omega e_1, \csc \theta \omega e_2,...,\csc \theta
\omega e_n\}$$   an adapted frame for slant
Riemannian maps.\\

We note that since the distribution $ker F_*$ is integrable it
follows that $\mathcal{T}_XY=\mathcal{T}_YX$ for $X, Y \in
\Gamma(ker F_*)$. Then the following Lemma can be obtained by using
Theorem 3.1.

\noindent{\bf Lemma 3.1.~}{\it Let $F$ be a slant Riemannian map
from a K\"{a}hler manifold to a Riemannian manifold. If $\omega$ is
parallel with respect to $\nabla$ on $ker F_*$, then}
\begin{equation}
\mathcal{T}_{\phi X}\phi X=-\cos^2 \theta \mathcal{T}_X X
\label{eq:3.7}
\end{equation}
{\it for $X \in \Gamma(ker F_*)$.}\\

In fact, proof of the above Lemma is exactly the same with the Lemma~3.3
given in \cite{Sahin1}.\\

We now give necessary and sufficient conditions for $F$ to be harmonic.\\

\noindent{\bf Theorem~3.2.~} {\it Let $F$ be a slant Riemannian map
from a K\"{a}hler manifold to a Riemannian manifold. Then $F$ is
harmonic if and only if}
\begin{equation}
\mathcal{T}_{\phi e_i}\phi e_i=-cos^2\, \theta\,\mathcal{T}_{e_i}e_i,
\label{eq:3.8}
\end{equation}
\begin{equation}
trace\mid_{\omega (ker F_*)}{^*F}_*(\mathcal{S}_{E_j}F_*(.))\in
\Gamma(\mu), \label{eq:3.9}
\end{equation}
{\it and}
\begin{equation}
trace\mid_{\mu}{^*F}_*(\mathcal{S}_{E_j}F_*(.))\in \Gamma(\omega(ker
F_*)), \label{eq:3.10}
\end{equation}
{\it where $\{e_1,\sec \theta \phi e_1, e_2,\sec \theta \phi e_2$
$,...,e_n,\sec \theta \phi e_n\}$ is an orthonormal frame for
$\Gamma(ker F_*)$ and $\{E_k\}$ is an orthonormal frame of
$\Gamma((range
F_*)^\perp)$..}\\

\noindent{\bf Proof.~}We choose a canonical orthonormal frame $e_1,
\sec \theta \phi e_1 ...,e_p, \sec \theta \phi e_p$, $\omega \csc
\theta e_{1}, ...,\omega \csc \theta e_{2p}$,
$\bar{e}_{1},...,\bar{e}_n$ such that $\{ e_1, \sec \theta \phi e_1
...,e_p, \sec \theta \phi e_p \}$ is an orthonormal basis of $ker
F_*$ and $\{ \bar{e}_{1},..,\bar{e}_{n}\}$ of $\mu$, where $\theta$
is the slant angle. Then $F$ is harmonic if and only if
\begin{eqnarray}
\sum^p_{i=1}(\nabla F_*)(e_i,e_i)&+&\sec^2\,\theta (\nabla F_*)(\phi
e_i,\phi e_i)+\csc^2\,\theta \sum^{2p}_{i=1}(\nabla F_*)(\omega
e_i,\omega e_i)\nonumber\\
&&+\sum^m_{j=1}(\nabla F_*)(\bar{e}_j,\bar{e}_j)=0.
\end{eqnarray}
By using (\ref{eq:2.1}) and (\ref{eq:2.7}) we have
\begin{equation}
\sum^p_{i=1}((\nabla F_*)(e_i,e_i)+\sec^2\,\theta (\nabla F_*)(\phi
e_i,\phi e_i)=-F_*(\mathcal{T}_{e_i}e_i+\sec^2\,\theta
\mathcal{T}_{\phi e_i}\phi e_i).\label{eq:3.11}
\end{equation}

On the other hand from Lemma 2.1, we know $\csc^2\,\theta
\sum^{2p}_{i=1}(\nabla F_*)(\omega e_i,\omega
e_i)+\sum^m_{j=1}(\nabla F_*)(\bar{e}_j,\bar{e}_j)\in \Gamma((range
F_*)^\perp)$. Thus we can write
\begin{eqnarray}
&&\csc^2\,\theta \sum^{2p}_{i=1}(\nabla F_*)(\omega e_i,\omega
e_i)+\sum^m_{j=1}(\nabla F_*)(\bar{e}_j,\bar{e}_j)= \csc^2\,\theta
\nonumber\\
&&\sum^{2p}_{i=1}\sum^{s}_{k=1}g_2((\nabla F_*)(\omega e_i,\omega
e_i),E_k)E_k\nonumber\\
&&+\sum^m_{j=1}\sum^{s}_{k=1}g_2((\nabla
F_*)(\bar{e}_j,\bar{e}_j),E_k)E_k\nonumber
\end{eqnarray}
where$\{E_k\}$ is an orthonormal basis of $\Gamma((range
F_*)^\perp)$. Then using (\ref{eq:2.12}) we have
\begin{eqnarray}
&&\csc^2\,\theta \sum^{2p}_{i=1}(\nabla F_*)(\omega e_i,\omega
e_i)+\sum^m_{j=1}(\nabla F_*)(\bar{e}_j,\bar{e}_j)= \csc^2\,\theta
\nonumber\\
&&\sum^{2p}_{i=1}\sum^{s}_{k=1}g_2(\mathcal{S}_{E_k}F_*(\omega
e_i),F_*(\omega
e_i))E_k\nonumber\\
&&+\sum^m_{j=1}\sum^{s}_{k=1}g_2(\mathcal{S}_{E_k}F_*(\bar{e}_j),F_*(\bar{e}_j))E_k\label{eq:3.12}
\end{eqnarray}
Then proof comes from the adjoint of $F_*$, (\ref{eq:3.11}) and
(\ref{eq:3.12}). \\

\noindent{\bf Example~5.~}Consider the slant Riemannian map given in
Example 4, then we have
$$(ker F_*)=Span\{Z_1=\frac{\partial}{\partial x_1},\,Z_2=\sin\,\alpha\,\frac{\partial}{\partial x_2}-\frac{\partial}{\partial x_3}+\cos\,\alpha\,\frac{\partial}{\partial x_4}\}$$
and
\begin{eqnarray}
(ker F_*)^\perp=Span\{Z_3&=&\frac{\sin\,\alpha}{\sqrt{2}}\,\frac{\partial}{\partial x_2}+\frac{1}{\sqrt{2}}\,\frac{\partial}{\partial x_3}+\frac{\cos\,\alpha}{\sqrt{2}}\,\frac{\partial}{\partial x_4},\nonumber\\
Z_4&=&-\cos\,\alpha\,\frac{\partial}{\partial x_2}+\sin\,\alpha\,\frac{\partial}{\partial x_4}\}.\nonumber
\end{eqnarray}
By direct computations, we have
$$JZ_1=-\frac{1}{2}Z_2+\frac{1}{\sqrt{2}}Z_3,\,JZ_2=Z_1+Z_4$$
which imply that
$$\phi Z_1=-\frac{1}{2}Z_2, \phi Z_2=Z_1.$$
Then it is easy to see that
$$\phi^2 Z_i=-\cos^2\frac{\pi}{4}Z_i=-\frac{1}{2}Z_i, i=1,2$$
which is the statement of Theorem~3.1. On the other hand, since $\mathcal{T}$ and $\mathcal{S}$ vanish for this slant Riemannian map, it satisfies the claim of Theorem~3.2.\\

By using (\ref{eq:2.7}) and (\ref{eq:3.3}), one can notice that the
equality (\ref{eq:3.8}) is satisfied in terms of
the tensor field $\omega$. More precisely, we have the following.\\

\noindent{\bf Lemma~3.2.~}{\it Let $F$ be a slant Riemannian map
from a K\"{a}hler manifold $(M_1,g_1,J)$ to a Riemannian manifold
$(M_2,g_2)$. If $\omega$ is parallel then (\ref{eq:3.8}) is
satisfied.}\\

\noindent{\bf Remark~1.~} We note that the equality (\ref{eq:3.7}) ( as a result of above lemma, parallel $\omega$) is enough for a slant submersion to be harmonic, however for a slant Riemannian map this case is not valid anymore.\\

We now investigate necessary and sufficient conditions for a slant
Riemannian map $F$ to be totally geodesic. We recall that a
differentiable map $F$ between Riemannian manifolds $(M_1,g_1)$ and
$(M_2,g_2)$ is called a totally geodesic map if
$(\nabla F_*)(X,Y)=0$ for all $X, Y \in \Gamma(TM_1)$. A geometric interpretation of a totally geodesic map is that it maps every geodesic in the total manifold into a geodesic in the base manifold in proportion to arc lengths.\\

\noindent{\bf Theorem~3.3.~}{\it Let $F$ be a slant Riemannian map
from a K\"{a}hler manifold $(M_1,g_1,J)$ to a Riemannian manifold
$(M_2,g_2)$. Then $F$ is totally geodesic if and only if}
$$
g_1(\mathcal{T}_U\omega V,\mathcal{B}X)=-g_2((\nabla F_*)(U,\omega
\phi V),F_*(X))+g_2((\nabla F_*)(U,\omega V),F_*(\mathcal{C}X))$$
$$
g_1(\mathcal{A}_X\omega U,\mathcal{B}Y)=g_2(\nabla^F_XF_*(\omega
\phi U),F_*(Y))-g_2(\nabla^F_XF_*(\omega U),F_*(\mathcal{C}Y))
$$
{\it and}
$$\nabla^F_XF_*(Y)+F_*(\mathcal{C}(\mathcal{A}_X\mathcal{B}Y+\mathcal{H}\nabla^1_X\mathcal{C}Y)+\omega(\mathcal{V}\nabla^1_X\mathcal{B}Y+\mathcal{A}_X\mathcal{C}Y))\in
\Gamma(range F_*)$$ {\it for $X, Y \in \Gamma((ker F_*)^\perp)$ and
$U,V \in \Gamma(ker F_*)$, where $\nabla^1$ is the Levi-Civita connection of $M_1$.}\\

\noindent{\bf Proof.~} From the decomposition of the total manifold
of a slant Riemannian map, it follows that $F$ is totally geodesic
if and only if $g_2((\nabla F_*)(U,V),F_*(X))=0$, $g_2((\nabla
F_*)(X,U),F_*(Y))=0$ and $(\nabla F_*)(X,Y)=0$ for $X, Y \in
\Gamma((ker F_*)^\perp)$ and $U,V \in \Gamma(ker F_*)$. First, since
$F$ is a Riemannian map, from (\ref{eq:2.1}) we obtain
$$g_2((\nabla F_*)(U,V),F_*(X))=-g_1(\nabla^1_UV,X).$$
Since $M_1$ is a K\"{a}hler manifold, using (\ref{eq:3.1}) and
(\ref{eq:3.2}) we have
\begin{eqnarray}
g_2((\nabla F_*)(U,V),F_*(X))&=&-\cos^2\,
\theta\,g_1(\nabla^1_UV,X)+g_1(\nabla^1_U\omega \phi V,X)\nonumber\\
&-&g_1(\nabla^1_U\omega V,\mathcal{B}X)-g_1(\nabla^1_U\omega
V,\mathcal{C}X).\nonumber
\end{eqnarray}
Taking into account that $F$ is a Riemannian map, using again
(\ref{eq:2.1}) and (\ref{eq:2.8}) we get

\begin{eqnarray}
g_2((\nabla F_*)(U,V),F_*(X))&=&\sec^2\,\theta
\{-g_1(\mathcal{T}_U\omega V,\mathcal{B}X)-g_2((\nabla F_*)(U,\omega \phi V),F_*(X))\nonumber\\
&+&g_2((\nabla F_*)(U,\omega V),F_*(\mathcal{C}X))\}.\label{eq:3.14}
\end{eqnarray}
In a similar way, we also have
\begin{eqnarray}
g_2((\nabla F_*)(X,U),F_*(Y))&=&\sec^2\,\theta
\{-g_1(\mathcal{A}_X\omega U,\mathcal{B}Y)-g_2(\nabla^F_XF_*(\omega U),F_*(\mathcal{C}Y))\nonumber\\
&+&g_2(\nabla^F_XF_*(\omega \phi U),F_*(Y))\}.\label{eq:3.15}
\end{eqnarray}
On the other hand, by using (\ref{eq:2.1}) and (\ref{eq:2.11}) we
derive
$$(\nabla F_*)(X,Y)=\nabla^F_XF_*(Y)+F_*(J\nabla^1_XJY)$$
for $X, Y \in \Gamma((ker F_*)^\perp)$. Then using (\ref{eq:3.1}),
(\ref{eq:3.2}) and (\ref{eq:2.7})-(\ref{eq:2.10}) we obtain
\begin{eqnarray}
(\nabla
F_*)(X,Y)&=&\nabla^F_XF_*(Y)+F_*(\mathcal{B}\mathcal{A}_X\mathcal{B}Y\nonumber\\
&+&\mathcal{C}\mathcal{A}_X\mathcal{B}Y+\phi
\mathcal{V}\nabla^1_X\mathcal{B}Y+\omega
\mathcal{V}\nabla^1_X\mathcal{B}Y\nonumber\\
&+&\mathcal{B}\mathcal{H}\nabla^1_X\mathcal{C}Y+\mathcal{C}\mathcal{H}\nabla^1_X\mathcal{C}Y\nonumber\\
&+&\phi\mathcal{A}_X\mathcal{C}Y+\omega\mathcal{A}_X\mathcal{C}Y).\nonumber
\end{eqnarray}
Since $$\mathcal{B}\mathcal{A}_X\mathcal{B}Y+\phi
\mathcal{V}\nabla^1_X\mathcal{B}Y+\mathcal{B}\mathcal{H}\nabla^1_X\mathcal{C}Y+\phi\mathcal{A}_X\mathcal{C}Y
\in \Gamma(ker F_*),$$
we have
\begin{eqnarray}
(\nabla
F_*)(X,Y)&=&\nabla^F_XF_*(Y)+F_*(\mathcal{C}\mathcal{A}_X\mathcal{B}Y\nonumber\\
&+&\omega \mathcal{V}\nabla^1_X\mathcal{B}Y+\mathcal{C}\mathcal{H}\nabla^1_X\mathcal{C}Y\nonumber\\
&+&\omega\mathcal{A}_X\mathcal{C}Y).\label{eq:3.16}
\end{eqnarray}
Then proof comes from (\ref{eq:3.14}), (\ref{eq:3.15}) and
(\ref{eq:3.16}).\\

\noindent{\bf Remark~2.~} Since $\mathcal{T}$, $\mathcal{A}$ and
$(\nabla F_*)$ vanish, Example 5 satisfies the conditions of Theorem
3.3.\\

\noindent{\bf Remark~3.~}We observe that the conditions for a slant
Riemannian map to be a totally geodesic are different from  the
conditions for a slant submersion to be totally geodesic, compare
Theorem 3.3 of the present paper with Theorem 3.5 of \cite{Sahin1}.
For a Riemannian submersion, the second fundamental form satisfies
$(\nabla F_*)(X,Y)=0$, $X, Y\in \Gamma((ker F_*)^\perp)$. However,
for a slant Riemannian map there is no guarantee that $(\nabla
F_*)(X,Y)=0$, $X, Y\in \Gamma((ker F_*)^\perp)$. From Lemma 2.1, we
only know that $(\nabla F_*)(X,Y)$ is $\Gamma((range F_*)^\perp)-$
valued. From the above reason it is necessary to use extra geometric
conditions to investigate the geometry of slant Riemannian maps.

\section*{4.~A decomposition theorem via slant Riemannian maps}
  \setcounter{equation}{0}
\renewcommand{\theequation}{4.\arabic{equation}}

In this section we are going to obtain necessary and sufficient
conditions for the total manifold of a slant Riemannian map to be a
locally product Riemannian manifold. Let $g$ be a Riemannian metric
tensor on the manifold $M = B \times F$ and assume that the
canonical foliations $D$ and $\bar{D}$ intersect perpendicularly
everywhere. Then from de Rham's theorem \cite{De}, we know that  $g$
is the metric tensor of a usual product Riemannian manifold if and
only if
$D$ and $\bar{D}$ are totally geodesic foliations. \\

\noindent{\bf Theorem~4.1.~}{\it Let $F$ be a slant Riemannian map
from a K\"{a}hler manifold $(M_1,g_1,J)$ to a Riemannian manifold
$(M_2,g_2)$. Then $(M_1,g_1)$ is a locally product Riemannian
manifold if and only if}
$$
g_1(\mathcal{T}_U\omega V,\mathcal{B}X)=-g_2((\nabla F_*)(U,\omega
\phi V),F_*(X))+g_2((\nabla F_*)(U,\omega V),F_*(\mathcal{C}X))$$

{\it and}
$$
g_2((\nabla F_*)(X,\mathcal{B}Y),F_*(\omega
U))=g_2(F_*(Y),\nabla^F_X F_*(\omega \phi
U))-g_2(F_*(\mathcal{C}Y),\nabla^F_X F_*(\omega U))
$$
{\it for $X, Y \in \Gamma((ker F_*)^\perp)$ and
$U,V \in \Gamma(ker F_*)$.}\\

\noindent{\bf Proof.~~} For $X, Y\in \Gamma((ker F_*)^\perp)$ and $U
\in \Gamma(ker F_*)$, from (\ref{eq:2.11}), (\ref{eq:3.1}),
(\ref{eq:3.2}) and Theorem 3.1, we have
\begin{eqnarray}
g_1(\nabla^1_XY,U)&=&-\cos^2\,
\theta\,g_1(Y,\nabla^1_XU)+g_1(Y,\nabla^1_X\omega \phi U)\nonumber\\
&-&g_1(\mathcal{B}Y,\nabla^1_X\omega
U)-g_1(\mathcal{C}Y,\nabla^1_X\omega U).\nonumber
\end{eqnarray}
Taking into account that $F$ is a Riemannian map and using
(\ref{eq:2.1}) we obtain
\begin{eqnarray}
g_1(\nabla^1_XY,U)&=&\sec^2\,\theta \{-g_2(F_*(Y),(\nabla
F_*)(X,\omega
\phi U))+g_2(F_*(Y),\nabla^F_XF_*(\omega \phi U))\nonumber\\
&-&g_2((\nabla F_*)(X,\mathcal{B}Y),F_*(\omega U))+g_2((\nabla
F_*)(X,\omega U),F_*(\mathcal{C}Y))\nonumber\\
&-&g_2(F_*(\mathcal{C}Y),\nabla^F_X F_*(\omega U))\}.\nonumber
\end{eqnarray}
Then Lemma 2.1 implies that
\begin{eqnarray}
g_1(\nabla^1_XY,U)&=&\sec^2\,\theta \{g_2(F_*(Y),\nabla^F_XF_*(\omega \phi U))-g_2((\nabla F_*)(X,\mathcal{B}Y),F_*(\omega U))\nonumber\\
&-&g_2(F_*(\mathcal{C}Y),\nabla^F_X F_*(\omega U)).\label{eq:4.1}
\end{eqnarray}
Thus proof follows from (\ref{eq:3.14}) and (\ref{eq:4.1}).\\

\noindent{\bf Acknowledgment.} The author is grateful to the referees for their valuable comments and suggestions. This paper is supported by The
Scientific and Technological Council of Turkey (TUBITAK) with number
(TBAG-109T125).

\end{document}